\documentclass[12pt,reqno]{article}

\usepackage{amsmath}

\usepackage{amsmath}
\usepackage{amsthm}
\usepackage{amscd}
\usepackage{amssymb}
\usepackage{graphicx}
\usepackage{epsf}
\theoremstyle{plain}
\newtheorem{theorem}{Theorem}[section]
\newtheorem{proposition}[theorem]{Proposition}
\newtheorem{lemma}[theorem]{Lemma}
\newtheorem{corollary}[theorem]{Corollary}

\newtheorem{definition}[theorem]{Definition}
\newtheorem{remark}[theorem]{Remark}
\newtheorem{example}[theorem]{Example}

\newtheorem{  -}[theorem]{  -}

\def\int{\text{int}}
\def\invlimit{\smash{\lim\limits_{\raise1pt\hbox{$\longleftarrow$}}}\vphanto
m{\big(}}
\def\inter{\hskip 1.5pt\raise4pt\hbox{$^\circ$}\kern -1.6ex}
\def\skel(#1,#2){#1^{(#2)}}
\def\hyp {\hbox {\rm {H \kern -2.8ex I}\kern 1.25ex}}
\def\reals {\hbox {\rm {R \kern -2.8ex I}\kern 1.15ex}}
\def\integers {\hbox {\rm { Z \kern -2.8ex Z}\kern 1.15ex}}
\def\naturals {\hbox {\rm {N \kern -2.8ex I}\kern 1.20ex}}
\def\rationals {\hbox {\rm { Q \kern -2.2ex l}\kern 1.15ex}}
\def\hyp {\hbox {\rm {H \kern -2.7ex I}\kern 1.25ex}}

\def\strutdepth{\dp\strutbox}
\def \ss{\strut\vadjust{\kern-\strutdepth \sss}}
\def \sss{\vtop to \strutdepth{
\baselineskip\strutdepth\vss\llap{$\diamondsuit\;\;$}\null}}

\def\strutdepth{\dp\strutbox}
\def \sst{\strut\vadjust{\kern-\strutdepth \ssss}}
\def \ssss{\vtop to \strutdepth{
\baselineskip\strutdepth\vss\llap{$\spadesuit\;\;$}\null}}

\begin{document}
\title{A finiteness result for Heegaard splittings}

\author{ Martin Lustig
\thanks  {supported by the Subvention  01/00495 for international
scientific cooperation of the R\'egion
Provence (France)}
\hskip9pt and Yoav Moriah\thanks{supported by The Fund for
Promoting Research at the Technion,  grant 100-127
and the Technion VRP fund, grant 100-127.}}

\date{}

\maketitle

\begin{abstract}

In this paper we show that for a given $3$-manifold and a given Heegaard splitting there are finitely
many preferred decomposing systems of $3g - 3$ disjoint essential disks. These are characterized by a
combinatorial criterion which is a slight strengthening of Casson-Gordon's rectangle condition. This is
in contrast to fact that in general there can exist infinitely many such systems of disks which satisfy just
the Casson-Gordon rectangle condition.

\quote{{\it Keywords:} Heegaard splittings, rectangle condition, double rectangle condition, pair of
pants decomposition

{\it 2000 AMS classification } 57N25 }

\end{abstract}

\vskip12pt

\section{Introduction}
\label{intro}

\vskip2pt

Every closed orientable 3-dimensional manifold $M$ admits a {\em Heegaard splitting}, i.e. a
decomposition into two handlebodies $H_{1}$ and $H_{2}$ which meet along their boundary.
This common boundary is called  a {\em Heegaard surface} in $M$ and is usually considered
only up to isotopy in $M$.

Heegaard splittings are a convenient way to define a $3$-manifold, but a priori it is difficult to get
structural information about the manifold from them. In the last fifteen years a lot of progress
was made in understanding the structural aspects of Heegaard splittings. A breakthrough was
achieved in the work of Casson and Gordon [CG] which ties Heegaard splittings to the existence
of incompressible surfaces. In particular, for non-Haken 3-manifolds strongly irreducible Heegaard
surfaces are now considered as suitable analogues of essential surfaces in the Haken case, thus
establishing them as an important tool in the study of these manifolds.

The main difficulty with Heegaard splittings is that a Heegaard splitting corresponds to a double coset
${\cal H} \phi {\cal H}$ of an element $\phi$ in the mapping class group ${\cal MCG}(\Sigma_{g})$ of 
a closed surface $\Sigma_{g}$ of genus $g \geq 2$, where  ${\cal H}$  is the subgroup of surface 
homeomorphisms which extend to a  handlebody $H$ via a properly chosen identifiction 
$\Sigma_{g} = \partial H$. This subgroup is not normal in ${\cal MCG}(\Sigma_{g})$, and it is not 
well understood at all.  The geometric analogue of this problem is the absence of a canonical \lq\lq coordinate 
system", that is a preferred choice of disks which define the handle structure in each of the two handlebodies 
of the splitting.

It is this problem that we wish to address. We choose a {\it complete decomposing system} ${\cal D}$,
of $3g - 3$, $g \geq 2$, disjoint non-parallel essential disks for each of the two handlebodies. These systems
${\cal D}_1 \subset H_1$ and ${\cal D}_2 \subset H_2$ decompose each of the handlebodies into
$2g - 2$ solid pairs of pants. Thus we obtain a {\it Heegaard diagram} for $M$, i.e. a finite set of
combinatorial data which determine $M$. There are infinitely many  such distinct complete decomposing
systems in each handlebody, so that the idea to recover characteristic data for $M$ from a Heegaard diagram
might seem hopeless. It is in this light that the following main result of this paper should be seen:

\vskip10pt

\noindent{\bf Theorem \ref{mainthm}.}
{\it For any closed orientable 3-manifold $M$ and any Heegaard splitting
$M = H_{1} \cup_{\partial H_{1} = \partial H_{2}} H_{2}$ there are only finitely many pairs of
complete decomposing systems ${\cal D}_{1} \subset H_{1}$ and ${\cal D}_{2} \subset H_{2}$
which satisfy the double rectangle condition.}

\vskip10pt

The double rectangle condition, defined precisely in Section \ref{doreco} below, is a slight
strengthening of the rectangle condition introduced by Casson and Gordon in [CG]. The statement that
Casson-Gordon's rectangle condition is generic, can be given a precise meaning using Thurston's measure
on the boundary of Teichm\"uller space. The question, whether the existence of complete decomposing
systems which satisfy the double rectangle condition is a generic property for Heegaard splittings, is
at present open (see Remark \ref{yettobewritten}).

As a corollary we obtain:
\vskip10pt

\begin{corollary}
\label{fmc} Let $M$ be an atoroidal closed $3$-manifold which admits a Heegaard splitting with two complete
decomposing systems that satisfy the double rectangle condition. Then the mapping class group of $M$ is  finite.

\end{corollary}

\begin{proof} It follows from a result of 	Jaco and Rubinstein [JR] that an atoroidal $3$-manifold has only
finitely many Heegaard splittings of any given genus. Any self-homeomorphism of $M$ must take  two complete
decomposing systems ${\cal D}_1,{\cal D}_2$ that satisfy the double rectangle condition to two other such
systems. and, by Theorem \ref{mainthm}, there are only finitely many of those. But every mapping class
which fixes ${\cal D}_1$ and ${\cal D}_2$ is easily seen to be trivial.

\end{proof}

\vskip0pt

\subsubsection*{Organization of the paper}

In Section \ref{doreco} we define the basic terminology and state our main result. We give a counterexample
to the conclusion of  Theorem \ref{mainthm} if the \lq\lq double rectangle condition" is replaced by the weaker 
\lq\lq Casson-Gordon rectangle condition". This shows that the  rectangle condition is not sufficient to characterize
a finite collection of  \lq\lq preferred" decomposing disk systems.

In Section \ref{disktypes} we investigate how the disks of a second complete decomposing system 
${\cal D}'_{1}$ in $H_{1}$ intersect the complementary components $B_{k}$ of the fixed
decomposing  system ${\cal D}_{1}$ in $H_{1}$ (these are solid pairs of pants). Any connected
component $\Delta' \subset {\cal D}'_{1}$ of this intersection is a disk which has as boundary an alternating
sequence of arcs from ${\cal D}_{1} \cap {\cal D}'_{1}$ and  from  ${\cal D}'_{1} \cap \partial H_{1}$.
The number of such arcs can be used as a measure of complexity for $\Delta'$. A priory there is no
bound on this complexity, which is one of the main reasons why homeomorphisms of 3-dimensional
handlebodies remain a mysterious and little understood topic.   In our context, however, one can exploit
the rectangle condition to get an upper bound on this complexity which depends on ${\cal D}_{1}$ and
${\cal D}_{2}$ only.  Even better, we show in Section \ref{disktypes} that, up to proper isotopy,
the disk $\Delta'$ must come from a finite collection which depends again only on ${\cal D}_{1}$
and ${\cal D}_{2}$.

In Section \ref{thickandthin} we investigate the complementary components of ${\cal D}'_{1}$ in each
solid pair $B_{k}$. They are called {\it parts}, and we distinguish {\it thin} and {\it thick} parts.
In the presence of the rectangle condition the possible nature and number of thick parts are both
determined by ${\cal D}_{1}$ and ${\cal D}_{2}$, while the number of thin parts depends
in an essential way also on ${\cal D}'_{1}$.

Bounding the number of the thin parts is the main problem in the proof of Theorem \ref{mainthm}
and is the only place where the double rectangle condition is used. This is accomplished in Section 5.

\vskip15pt

\begin{remark}
\label{similar} \rm The intersection pattern induced by the disks from ${\cal D}'_{1}$ on every solid pair
of pants $B_{k}$ is strongly reminiscent of the intersection pattern on a $3$-simplex given by a surface
$S$ in normal position, which is cut by $S$ into a bounded number of {\it thick} blocks and an arbitrary
number of {\it thin} pieces that occur in \lq\lq parallel stacks" (compare e.g. [Sch]). One important difference
is that normal surface theory is done for closed  surfaces, while we work with disks in handlebodies.

\end{remark}

\noindent {\bf Acknowledgments:} We would like to thank Hyam Rubinstein for pointing out
Corollary \ref{fmc} to us, and Saul Schleimer for many inspiring discussions. We also thank
the Universit\'e d'Aix-Marseille III and the Technion, where most of this work was done. Finally,
special thanks go to Caf\'e Parisien in Marseille for its hospitality.

\vskip25pt

\section{The double rectangle condition}
\label{doreco}

\vskip10pt

Let $M$ be a closed $3$-dimensional manifold, and let $\Sigma \subset  M$ be a closed orientable
{\it Heegaard surface} of genus $g \geq 2$ cutting $M$  into two handlebodies $H_{1}$ and $H_{2}$.

Let ${\cal D}_{1} \subset H_1$  and ${\cal D}_{2} \subset H_2$ be two complete decomposing disk systems,
i.e. each handlebody is decomposed by the disk system into a union of solid pairs of pants. We will always
assume that $\partial {\cal D}_{1}$ and $\partial {\cal D}_{2}$ have only essential intersections, that is,
they intersect in  transverse intersection points, and one can not decrease their number by a proper isotopy of
${\cal D}_{1}$ in $H_{1}$ or of ${\cal D}_{2}$ in $H_{2}$.

A {\it wave }  $\omega \subset \Sigma$ with respect to ${\cal D}_{1} $ is an arc that meets
${\cal D}_{1}$  only in its boundary points $\partial \omega$, which lie on the same component
$\partial D_j \subset \partial {\cal D}_1$, such that in $H_{1}$ the arc $\omega$ is  isotopic
relative endpoints to a subarc of $\partial D_j$, but not in $\Sigma$. Similarly we define waves
for ${\cal D}_2 $.

We say that the closure of a connected component of  $\Sigma - (\partial{\cal D}_{1} \cup \partial {\cal D}_{2})$
is a {\it rectangle} $R$ if it is homeomorphic to a disk, whose boundary $\partial R$ is a concatenation
of precisely four arcs, two of which are subarcs on curves in $\partial {\cal D}_{1}$ and the other two
are subarcs of curves  in $\partial {\cal D}_{2}$. It is possible that two of the curves from one system
belong to the same component, and even that two opposite \lq\lq boundary vertices" of the rectangle are
itentified.

An {\it adjacent pair of curves} in $\partial {\cal D}_{1}$ (similarly in ${\cal D}_{2}$) consists
of  two  curves which can be joined by an essential arc in $\Sigma - \partial {\cal D}_{1}$ which
does not meet other curves from $\partial {\cal D}_{1}$, and which is  not a wave.  Such an
arc lies in one of the pair of pants of the decomposition defined by $\partial {\cal D}_{1}$,
and is unique up to isotopy in this pair of pants, so that we usually suppress its mentioning
and only note the two curves in $\partial {\cal D}_{1}$. Similarly, an {\em adjacent triple of curves}
in $\partial {\cal D}_{1}$  consist of three curves which can be connected by an arc that intersects 
the middle curve transversely, and the resulting two subarcs define two adjacent pairs of curves. Note that
the above two definitions include the situation where the inclusion of the pair of pants into the surface
$\Sigma$ identifies two of its boundary curves. The same definitions hold for ${\cal D}_2 \subset H_2$.

\vskip10pt

A. Casson and C. Gordon have introduced the following [CG]:

\vskip10pt

\begin {definition}\rm
\label {Rectangle condition}
The complete decomposing systems ${\cal D}_{1} \subset H_{1}$ and ${\cal D}_{2} \subset H_{2}$
satisfy the {\it rectangle condition} if every pair of adjacent curves in $\partial {\cal D}_{1}$  and any
pair of adjacent curves in $\partial {\cal D}_{2}$  form at least one rectangle which is contained in the
intersection of the respective pairs of pants.

\end{definition}

\vskip0pt

The importance of this notion comes from Casson-Gordon's observation that a Heegaard splitting
$M = H_{1} \cup_{\Sigma} H_{2}$ which satisfies the rectangle condition is strongly irreducible:
Indeed, any essential  disk $D \subset H_{1}$ must either be parallel to a curve of ${\cal D}_{1}$
or contain a wave with respect to ${\cal  D}_{1}$.  In both cases there exist two adjacent curves of
${\cal D}_{1}$ such that $D$ intersects all rectangles formed by these two curves with any adjacent
pair of curves from  $\partial {\cal D}_{2}$. As the analogue is true for any essential disk 
$E \subset H_{2}$, it follows from the rectangle condition that $D$ and $E$ must intersect in one of the
rectangles, so that the Heegaard  splitting is strongly irreducible. In particular all waves with respect to
${\cal D}_{1}$ must intersect all waves with respect to ${\cal D}_{2}$.

\vskip8pt

The same idea is used in the proof of the next lemma.

\vskip10pt

\begin{lemma}
\label{nowaves}
(a) If  ${\cal D}_{1}$ and ${\cal D}_{2}$ satisfy the rectangle condition, then for any disk $D \subset H_{1}$
the boundary curve  $\partial D \subset \Sigma$ does not contain a wave with respect to ${\cal D}_{2} \subset H_{2}$.

\noindent (b) Every wave on $\Sigma$ with respect to ${\cal D}_1$ intersects every curve which bounds a disk in
$H_{2}$ at least once.

\end{lemma}

\begin{proof}
 (a) As ${\cal D}_{1}$ is a complete decomposing system of $H_{1}$, the curve $\partial D$
must either be parallel to one of the  $\partial D_{i}$, or it contains a wave with respect to
$\partial {\cal D}_{1}$.   In both cases there exist two adjacent curves of ${\cal  D}_{1}$ such
that $\partial D$  intersects all rectangles formed by these two  curves with any adjacent pair
of curves from  $\partial {\cal D}_{2}$.  If  $\partial D$ also contains a wave with respect to
$\partial {\cal D}_{2}$, then there exist two adjacent curves of ${\cal D}_{2}$ with the same property.
Hence we could deduce from the rectangle condition at least one self-intersection of $\partial D$
in one of the rectangles.

\noindent
(b)  The claim follows exactly from the same arguments.

\end{proof}

\vskip10pt

\begin{remark}
\label{nonfinite}\rm
It is possible that a given Heegaard splitting possesses infinitely many non-isotopic decomposing disk
systems ${\cal D}_1$ and $ {\cal D}_2$ all satisfying the rectangle condition.  An example will be given
at the end of this section.

\end{remark}

\vskip7pt

To get the desired finiteness result  Theorem \ref{mainthm}, we have to strengthen the rectangle
condition slightly:  We call the union of two rectangles which have a side in common, a {\it double rectangle}.
Thus the boundary of a double rectangle formed by ${\cal D}_1$ and ${\cal D}_2 $ consists of two 
subarcs from an adjacent pair of curves of,  say, $\partial {\cal D}_1$, and of two subarcs from the two 
outer curves of an adjacent triple of curves of $\partial {\cal D}_2$.

\vskip10pt

\begin{definition}
\label{drc}  \rm
The decomposing disk systems ${\cal D}_1$ and $ {\cal D}_2$  satisfy the {\it double rectangle condition}
if every pair of adjacent curves from $\partial {\cal D}_1$ forms, with every adjacent triple from
$\partial {\cal D}_2$, a double rectangle,  and vice versa.

\end{definition}

\vskip8pt

Note that, of course, the double rectangle condition implies the rectangle condition.

\vskip8pt

\begin{lemma}
\label{lemma1}
If ${\cal D}_1$ and ${\cal D}_2$ satisfy the double rectangle condition,
then every essential disk
$D \subset H_{1}$ intersects every triple in ${\cal D}_2$, and vice versa.

\end{lemma}

\begin{proof}
If $D$ belongs (perhaps after a proper
isotopy in $H_{1}$) to ${\cal D}_1$,
then the claim is obviously
true. Otherwise, the curve $\partial D$
has a wave with respect to ${\cal
D}_1$.  This implies that  there
is at least one adjacent pair of curves in some pair of pants in ${\cal
D}_1$  which is separated by this wave.
Since ${\cal D}_1$ and ${\cal D}_2$ satisfies the double rectangle
condition, the adjacent pair and  hence
the curve $\partial D$ must intersect any adjacent triple of curves from
${\cal D}_2$.

\end{proof}

\vskip10pt

It follows that on an adjacent pair of pants we have the following
intersection pattern as in Fig. 1.

\begin{figure}[htb]
\vbox{{\epsfysize240pt\epsfbox{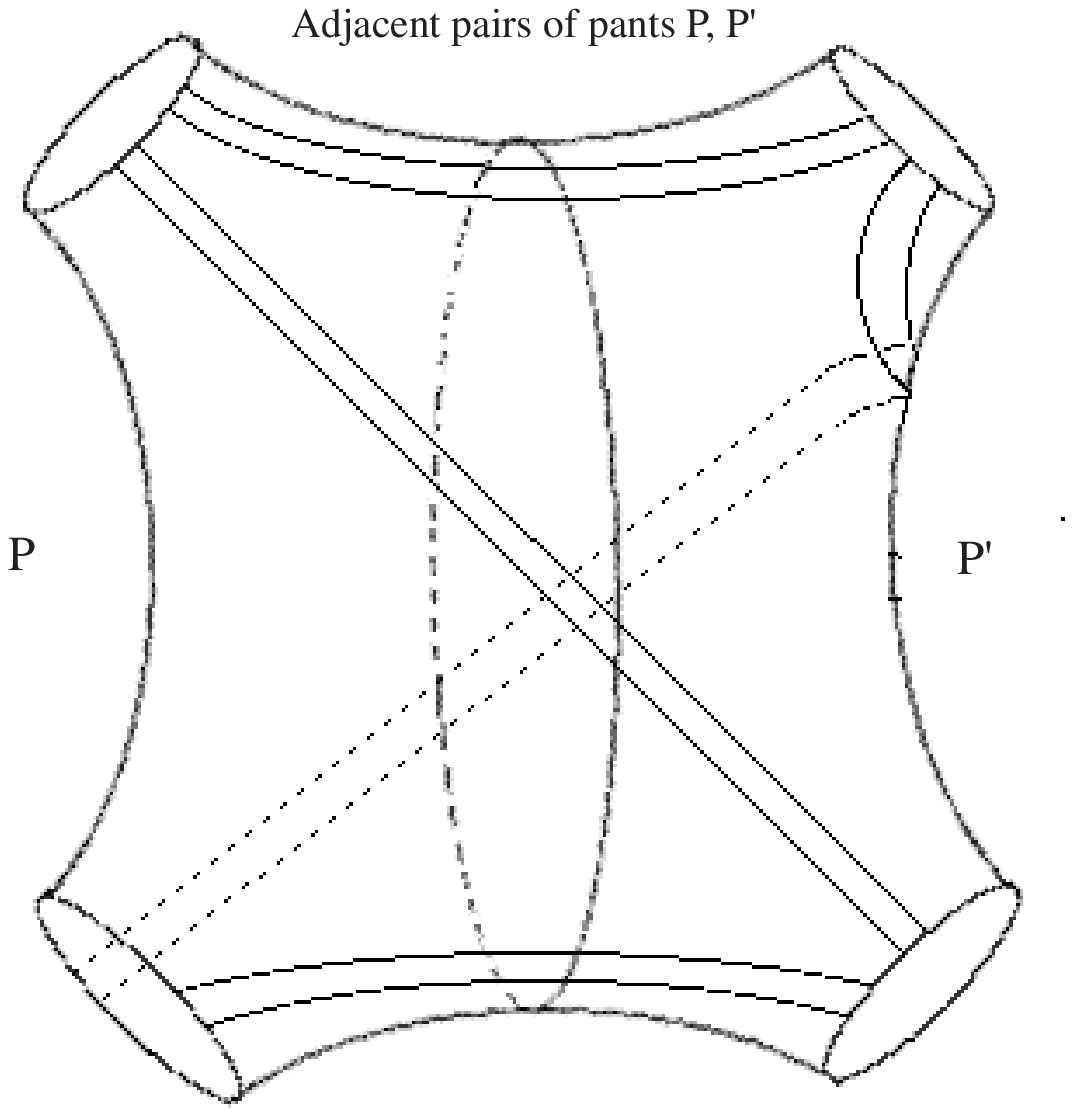}}}
\centerline {Figure 1.}
\end{figure}

\vskip10pt

\noindent We can now state the main result of this paper:

\vskip10pt

\begin {theorem}
\label{mainthm}
For any closed orientable 3-manifold $M$ and any Heegaard splitting $M = H_{1} \cup_{\Sigma} H_{2}$
there are only finitely many pairs of complete decomposing systems ${\cal D}_{1} \subset H_{1}$
and ${\cal D}_{2} \subset H_{2}$ which satisfy the double rectangle condition.

\end{theorem}

\vskip6pt

We finish this section with a counterexample to the analogue of this result, if one replaces
the double rectangle condition by the simple rectangle condition:

\begin{example}
\label{example1}
\rm Consider the genus two Heegaard diagram obtained from Fig. 2 by making the following
identifications:  $D_1 \equiv D'_1 \,,\,\, x \equiv x' \,,\,\,  y \equiv y'$ and
$D_2 \equiv D'_2 \,,\,\,  w \equiv w' \,,\,\,  z \equiv z'$.

Let $H_1$ be the genus two handlebody obtained by these identifications from Fig. 2, and let $H_2$ be an
identical copy of $H_{1}$. Let $M = H_1 \cup_{\delta^{t}} H_2$ , where $t$ is some sufficiently
large integer, and $\delta^{t}$ is the ${t}$-fold Dehn twist along the curve $\delta \subset \partial H_1$.
Let ${ \cal D }_1$ be the complete decomposing system given by the disks $\{D_1,D_2,D_3\}$  and
${ \cal D }_2$ be the identical system in $H_2$. Note that our choice of the Dehn twist exponent ensures that
the two systems  ${ \cal D }_1$ and  ${ \cal D }_2$ satisfy the rectangle condition.

Now consider the annulus $A \subset H_1$ as in Fig. 2 and change the system ${\cal D }_1$ to
a system ${ \cal D }^n_1$ by twisting $n$ times along $A$. It is immediate to see that all systems
${ \cal D }^n_1$ together with  the system ${ \cal D }_2$ satisfy the rectangle condition for all
$n \in \integers$.

\end{example}

\vbox{{\epsfysize190pt\epsfbox{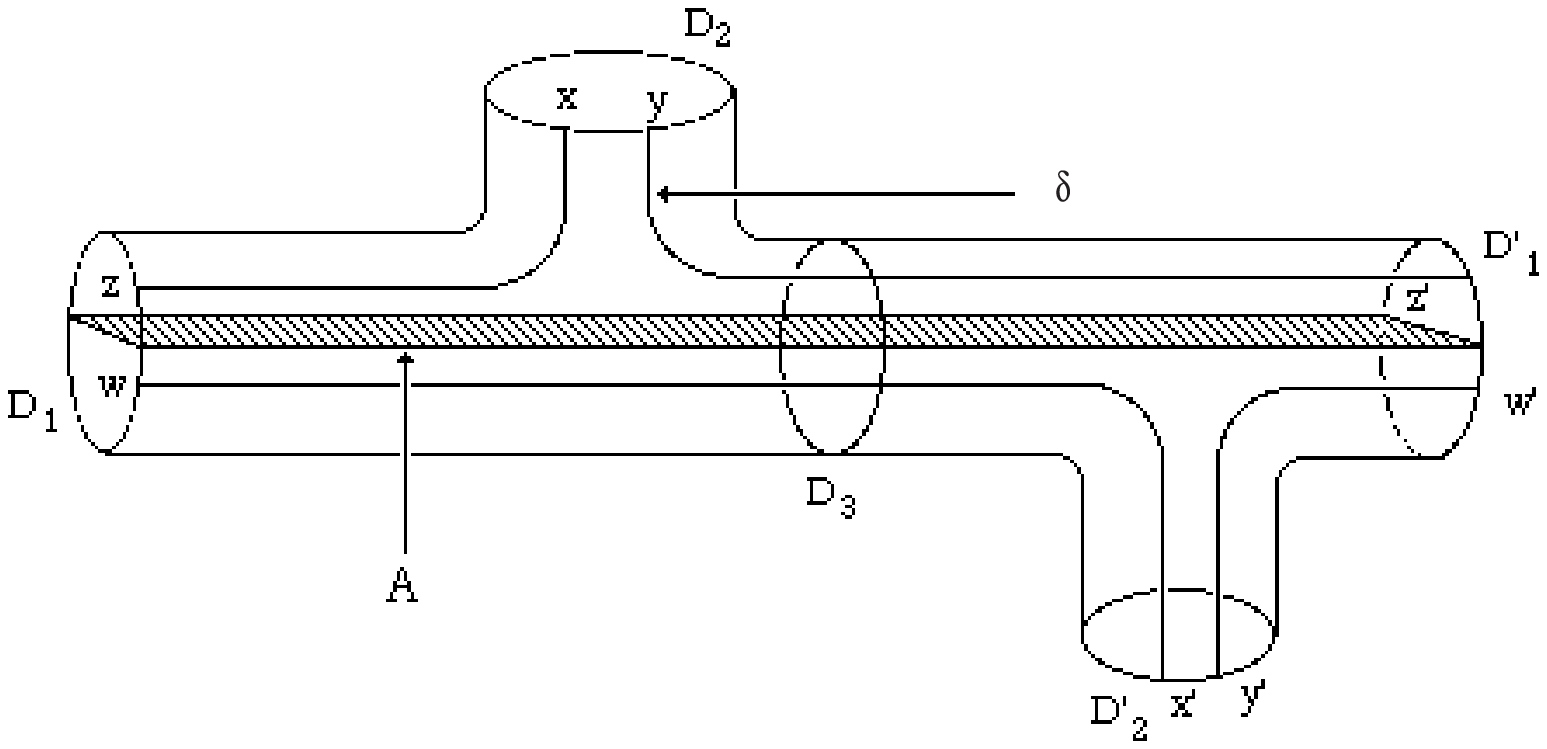}}}
\centerline {Figure 2.}

\vskip20pt

\section{Finiteness of disk types}
\label{disktypes}

\vskip10pt

We now concentrate on the handlebody $H_{1}$ which contains two complete decomposing disk systems
${\cal D}_{1}$ and ${\cal D}'_{1}$.  We think of ${\cal D}_{1}$ as being the fixed reference system, and
of ${\cal D}'_{1}$ as an alternative candidate: The goal of the paper is to show that, under the right conditions,
there are only finitely many such ${\cal D}'_{1}$.

\vskip6pt

In order to simplify the terminology we define:

\vskip6pt

\begin{definition}
\label{uniformly} \rm  We say that a constant defined by means of ${\cal D}'_{1}$ is {\it uniformly bounded}
if it depends only on the fixed pair of  decomposing systems of disks ${\cal D}_1 \subset H_1$ and
${\cal D}_2 \subset H_2$.

\end{definition}

\vskip0pt

As handlebodies are irreducible, we can assume that (after a suitable isotopy) ${\cal D}_{1}$ and
${\cal D}'_{1}$ are {\em tight}: They intersect only in arcs which terminate in essential
intersection points of their boundary curves. Thus each disk of ${\cal D}'_{1}$ is cut by
${\cal D}_{1}$ into {\em  disk pieces} which have as boundary an alternating sequence of
{\em intersection arcs} from ${\cal D}_{1} \cap{\cal D}'_{1}$ and {\em connecting arcs} from
$\partial {\cal D}'_{1} \subset \partial H_{1}$.

Every connecting arc is contained in a single pair of pants from the decomposition of $\partial H_{1}$ with
respect to ${\cal D}_{1}$, and  it can not be boundary parallel on this pair of pants: This follows from our
assumption that ${\cal D}_{1}$ and  ${\cal D}'_{1}$ are tight. For intersection arcs we prove the weaker
fact that  they can not be boundary parallel on ${\cal D}_{1} - \partial {\cal D}_{2}$:

\vskip10pt

\begin{lemma}
\label{lemma2}
Let ${\cal D}_1, {\cal D}'_1 \subset H_{1}$ and  ${\cal D}_2, {\cal D}'_2 \subset H_{2}$
be complete decomposing systems, and assume that the pair ${\cal D}'_1, {\cal D}'_2$ satisfies
the rectangle condition. Then for any disk $D_{i} \subset {\cal D}_1$ every intersection arc
$\alpha \subset D_{i}$ has its endpoints in two distinct connected components of
$\partial D_{i}- \partial {\cal D}_2$.

\end{lemma}

\begin{proof}
It suffices to consider an intersection arc $\alpha$ which is  contained in the boundary of an
outermost subdisk $\Delta$ of $D_{i} \in {\cal D}_1$. Every such $\Delta$ contains in its
boundary an arc  $\omega = \partial \Delta - \inter \alpha$ $\subset \partial D_{i}$.
As $\Delta$ is outermost, $\omega$ meets ${\cal D}'_{1}$ only in its boundary points, and
hence is a wave on $D_{i} $  with respect to ${\cal D}'_1$.

We can apply Lemma \ref{nowaves} (b) to ${\cal D}'_1$ and ${\cal D}'_2$ to conclude that $\omega$
must meet every curve of $\partial {\cal D}_2$ (see Fig. 3 below ).

\end{proof}

\vskip15pt

We now use the disks from ${\cal D}_{2}$ to group the intersection and the connecting arcs,
defined  above, into equivalence classes:  Given a disk $ D_{i} \subset {\cal D}_1$, two
intersection arcs $\alpha, \alpha ' \subset D_{i} \cap {\cal D}'_1$ will be called {\it parallel} if
the pair  $(\alpha,\partial \alpha)$ is isotopic to the pair $(\alpha',\partial  \alpha' )$ in
$(D_i, \partial D_i - \partial {\cal D}_2 )$. Similarly, two connecting arcs $\beta, \beta '$
will be call {\it parallel} if the pair $(\beta,\partial \beta)$ is isotopic to the pair $(\beta',\partial \beta' )$
in  $(\partial H_{1}, \partial {\cal D}_{1} - \partial {\cal D}_2)$. Such an isotopy class of parallel arcs
will be called the  {\it arc type} of an intersection arc or of a connecting arc.

\vbox{{\epsfysize200pt\epsfbox{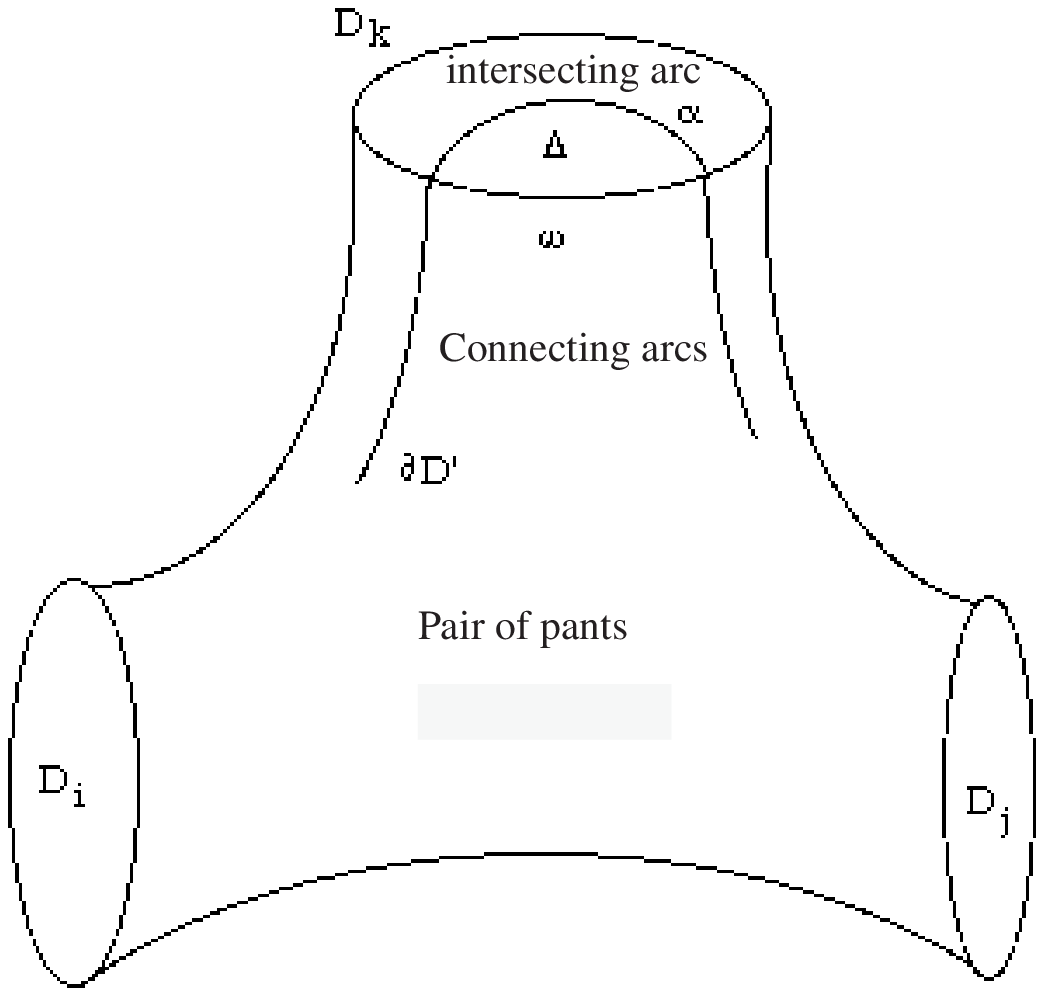}}}
\centerline {Figure 3.}

\vskip20pt

It follows from Lemma \ref{lemma2} and from the stronger fact for connecting arcs, stated in the
paragraph just before Lemma \ref{lemma2}, that two arcs $\alpha$ and $\alpha'$ which belong to
the same arc type are indeed parallel: They span a band (in $\partial H_{1}$ or in ${\cal D}_{1}$)
where the ``long" sides are given by $\alpha$ and $\alpha'$, while the ``short" sides are arcs from
$\partial {\cal D}_{1} - \partial {\cal D}_{2}$.

\vskip15pt

\begin{lemma}
\label{mularcs}
Let ${\cal D}_1, {\cal D}'_1 \subset H_{1}$ and ${\cal D}_2, {\cal D}'_2 \subset H_{2}$
be complete decomposing systems, and assume that the pair ${\cal D}'_1, {\cal D}'_2$ satisfies
the rectangle condition. Then  the number of intersection arc types on $\partial {\cal D}'_1$ with
respect to ${\cal D}_1$ is uniformly bounded above.

\end{lemma}

\begin{proof}

The system ${\cal D}_2$ determines the number of points of $\partial{\cal D}_2$ on each of the
$3g - 3$ components of $\partial{\cal D}_1$. Hence it determines their complementary components on
$\partial{\cal D}_1$. Thus there are finitely many relative isotopy classes of arcs (in $\partial H_{1}$
or in ${\cal D}_1$)  connecting them.

\end{proof}

\vskip8pt

\begin{lemma}
\label{mulconarcs}
Let ${\cal D}_1, {\cal D}'_1 \subset H_{1}$ and ${\cal D}_2, {\cal D}'_2 \subset H_{2}$
be complete decomposing systems, and assume that the pair ${\cal D}'_1, {\cal D}'_2$ satisfies
the rectangle condition. Then the number of connecting arc types on $\partial {\cal D}'_1$ with
respect to ${\cal D}_1$ is uniformly bounded above.

\end{lemma}

\begin{proof}
Every connecting arc $\alpha$ is contained in a single pair of pants $P$ from the decomposition of
$\partial H_{1}$ with respect to ${\cal D}_{1}$. Hence its isotopy class relative endpoints is essentially
determined by the choice  of the  boundary curves from $\partial P \subset \partial {\cal D}_1$ which
contain the endpoints of $\alpha$. More precisely, up to relative isotopy these arcs are determined by
the intervals on such a boundary curve which in turn are determined by the intersections with the
system ${\cal D}_{2}$, up to possible twists around these boundary curves. Thus we need to show
that there are only finitely many choices for the number of such twists:

As the connecting arcs are disjoint among themselves, if one of them spirals around a boundary component
$\partial D_i$ of $P$,  then so do all of those connecting arcs which have an endpoint on $\partial D_i$.
This spiraling is \lq\lq controlled" by the arcs from $\partial {\cal D}_2$ in $P$: By Lemma \ref{nowaves} (b)
for each $D_i$ from  ${\cal D}_1$ there must be at least one arc $\beta$ from $P \cap \partial {\cal D}_2$
which intersects $\partial D_i$.

We note that somewhere on  $\partial D_i$ there must be  a wave with respect to $\partial {\cal D}'_1$:
This wave is given by two adjacent intersection points on $\partial D_i$ with two connecting arcs
$\alpha_1, \alpha_2$  that lie on the same curve $\partial D'_j \subset \partial {\cal D}'_1$, such that, when
running once around $\partial D'_j$, the arcs $\alpha_1, \alpha_2$ are traversed in opposite directions (see Fig. 4).

\vskip25pt

\vbox{{\epsfysize140pt\epsfbox{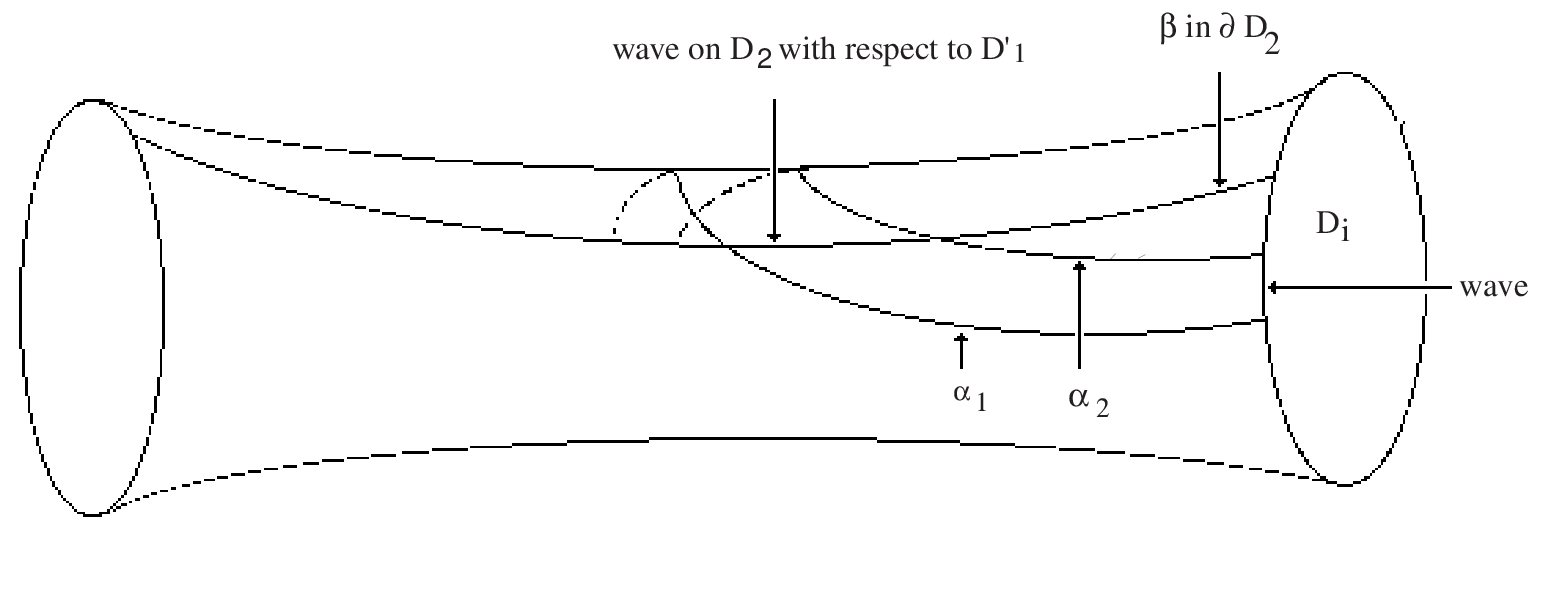}}}
\centerline{Figure 4.}

\vskip20pt

Now assume that $\alpha_1$ and $\alpha_2$ spiral around $\partial D_i$ for some time, in a parallel fashion,
thus intersecting the above  arc $\beta$ at least once. But then the band spanned by the spiraling arcs $\alpha$
and $\alpha'$ intersects $\beta$ in a wave on $\beta \subset \partial D_{k} \subset \partial{\cal D}_2$ with respect
to  $\partial {\cal D}'_1$. Since the disk $D_k$ belongs to ${\cal D}'_2$ or has some wave with  respect to
${\cal D}'_2$ this would contradict Lemma \ref{nowaves} (a). Hence $\alpha_1$ and $\alpha_2$ can not  spiral
around $\partial D_i$, and hence there are only finitely many connecting arc types on any pair of pants $P$
which are determined only by ${\cal D}_1$ and ${\cal D}_2$.

\end{proof}

\vskip10pt

We call the components of $H_1$, when cut along ${\cal D}_1$ (or  ${\cal D}'_1$), {\em solid pairs of pants}
and denote them by $B_k $ (or $B'_k $ resp.), for $k = 1,\dots 2g - 2$.  Denote by ${\cal B}$ (or ${\cal B}'$ resp.)
the  collection of these solid  pairs of pants. We defined above a disk piece to be a connected component of some 
${\cal D}'_1 \cap B_k$. Define a  {\it disk type} to be a class of disks pieces whose boundaries are composed of
intersection arcs and connecting arcs which are parallel pairwise. It follows from the previous discussion that disk pieces
which belong to the same disk type lie in one of the $B_{k}$ as a parallel stack, that is, homeomorphic to horizontal disks
in $D^2 \times \reals$.

A priori a disk piece can have in its boundary distinct connecting arcs or intersection arcs that belong to
the same arc type. However, this turns out to be impossible, if the rectangle condition is  imposed:

\vskip12pt

\begin{lemma}
\label{lemma3}
Let ${\cal D}_1, {\cal D}'_1 \subset H_{1}$ and ${\cal D}_2, {\cal D}'_2 \subset H_{2}$
be complete decomposing systems, and assume that the pair ${\cal D}'_1, {\cal D}'_2$ satisfies
the rectangle condition. Then any intersection arc type or connecting arc type can occur in the
boundary of a given disk piece at most once.

\end{lemma}

\begin{proof}
Given a disk piece $\Delta' \subset {\cal D}_1'$, orient its boundary $\partial \Delta'$ and assume that some
arc type appears more than once in $\partial \Delta'$.  Hence there are two distinct arcs $\alpha_1, \alpha_2$
in  $\partial \Delta'$ which belong to the same arc type.

Let $B_{k}$ be the solid pair of pants that contains $\Delta'$. Note that $\partial B_k$ is a $2$-sphere and
$\partial \Delta'$  is a simple closed curve on this sphere. Hence, if  the orientation induced on $\alpha_1$ and
$\alpha_2$ by the choice of orientation on $\partial \Delta'$ induces on them the same orientation as parallel
intersection or connecting arcs,  then there must be a third arc $\alpha_3$ in $\partial \Delta'$ of the same arc
type, such that  $\alpha_3$ runs  between $\alpha_1$ and $\alpha_2$, but with the opposite orientation: Otherwise
$\partial \Delta'$ would either not be simple or not be a  closed curve.

Hence we can assume by a standard innermost argument that $\alpha_1$ and $\alpha_2$ are adjacent
arcs in  the same arc type,  and that $\partial  \Delta$ traverses them in opposite directions.
Let $\partial D_i \subset \partial {\cal D}_1$ be the curve which contains an endpoint of this arc type i.e.,
$D_i \subset {\cal D}_1$ is one of the three boundary disks of $B_k$. Furthermore let $\beta$ be the
subarc on $\partial D_i$ which joins  the endpoints of $\alpha_1$ and $\alpha_2$. Since the two arcs
are adjacent in the arc type, and are traversed by $\partial \Delta'$ in opposite directions, it follows that
$\beta$ is a wave on $\partial D_i \subset \partial {\cal D}_1$  with respect to ${\cal D}'_1$. In particular,
$\beta$ does not meet ${\cal D}_2$ in it's interior. As we assume tha ${\cal D}'_1$ and  ${\cal D}'_2$
satisfy the rectangle condition, this contradicts Lemma \ref{nowaves} (b).

\end{proof}

\vskip10pt

\begin {proposition}
\label{finitetype}
Let ${\cal D}_1, {\cal D}'_1 \subset H_{1}$ and ${\cal D}_2, {\cal D}'_2 \subset H_{2}$
be complete decomposing systems, and assume that the pair ${\cal D}'_1, {\cal D}'_2$ satisfies
the rectangle condition. Then there is a finite set of disk types in any of the solid pair of pants $B_{k}$
from $H_{1} - {\cal D}_1$, such that any of the disk pieces of ${\cal D}'_1 - {\cal D}_1$ belongs to
one  of the disk types in the above finite set. Furthermore, the number of disk types in this finite set is
uniformely bounded above.

\end{proposition}

\vskip10pt

\begin{proof}
We can apply Lemma \ref{mularcs}, Lemma \ref{mulconarcs} and Lemma \ref{lemma3} to conclude
that ${\cal D}_1$ and ${\cal D}_2$ determine a finite set of intersecting arc types, and a finite set of connecting
arc types, which can possibly appear in the boundary of a disk type $\Delta$. Furthermore, each of those appears
in the boundary of  $\Delta$ at most once.  Hence there are only finitely many possible disk types for $\Delta$,
and they are dependent only on ${\cal D}_1$ and ${\cal D}_2$.

\end{proof}

\vskip10pt

\begin{remark} \rm Note that in all of
	Lemmas \ref{lemma2} - \ref{lemma3} and Proposition \ref{finitetype} we require that only ${\cal D}'_1$ 
and ${\cal D}'_2$ satisfy the rectangle condition,  but not necessarily ${\cal D}_1$ and ${\cal D}_2$.

\end{remark}

\vskip20pt

\section{Thick and thin regions}
\label{thickandthin}

\vskip15pt

In the last section we considered the solid pairs of pants $B_{k} \in {\cal B}$ obtained from cutting the
handlebody $H_1$ along the complete decomposing disk system ${\cal D}_1$. In this section we change
our point of view and consider the solid pairs of pants $B'_j$, obtained from cutting $H_1$ along the
disk system ${\cal D}'_1$. The collection of these solid pairs of pants will be called ${\cal B}'$. The
connected components of the intersection $B_k \cap B'_j$ of any $B_k \in {\cal B}_1$  with any
$B'_l \in {\cal B}'_1$ are called {\em parts}, and we distinguish  two kinds of them:

\vskip10pt

\begin{definition}\rm

A connected component of  $B_k \cap B'_l$ is called a {\it thin part} if its intersection with  ${\cal D}'_1$
consists of two disk pieces which belong to the same disk type in $B_i \in {\cal B}$.  Otherwise the
connected component is called a {\it thick part}.

\end{definition}

\vskip10pt

In any solid pair of pants $B_{k}$ a {\em stack} is a maximal collection of thin parts. The boundary of the
stack is composed of disk pieces from ${\cal D}'_1$ all belonging to the same disk type. Notice that the
complementary components  in $B_{k}$ of the union of all stacks are precisely the thick parts.

\vskip10pt

We now want to group together the parts in one solid pair of pants $B'_{l}$ into larger units, called {\em regions}:

\vskip10pt

\begin{definition} \rm
\label{regions}
For each $B'_l$, a {\em thick region} is a maximal union of thick parts in $B'_l$ which is connected.
A region is {\em thick peripheral} if it is disjoint from at least one of the three boundary disks of  $B'_i$
from the system ${\cal D}'_1$  (see Fig. 7). The region is called {\it central} if all three boundary disks
are met (see Fig. 5). A {\it thin region} is a maximal connected union of thin parts contained in $B'_l$
(see Fig. 6).  Finally, the {\em volume} of any region is the number of parts contained in that region.

\end{definition}

\vbox{{\epsfysize180pt\epsfbox{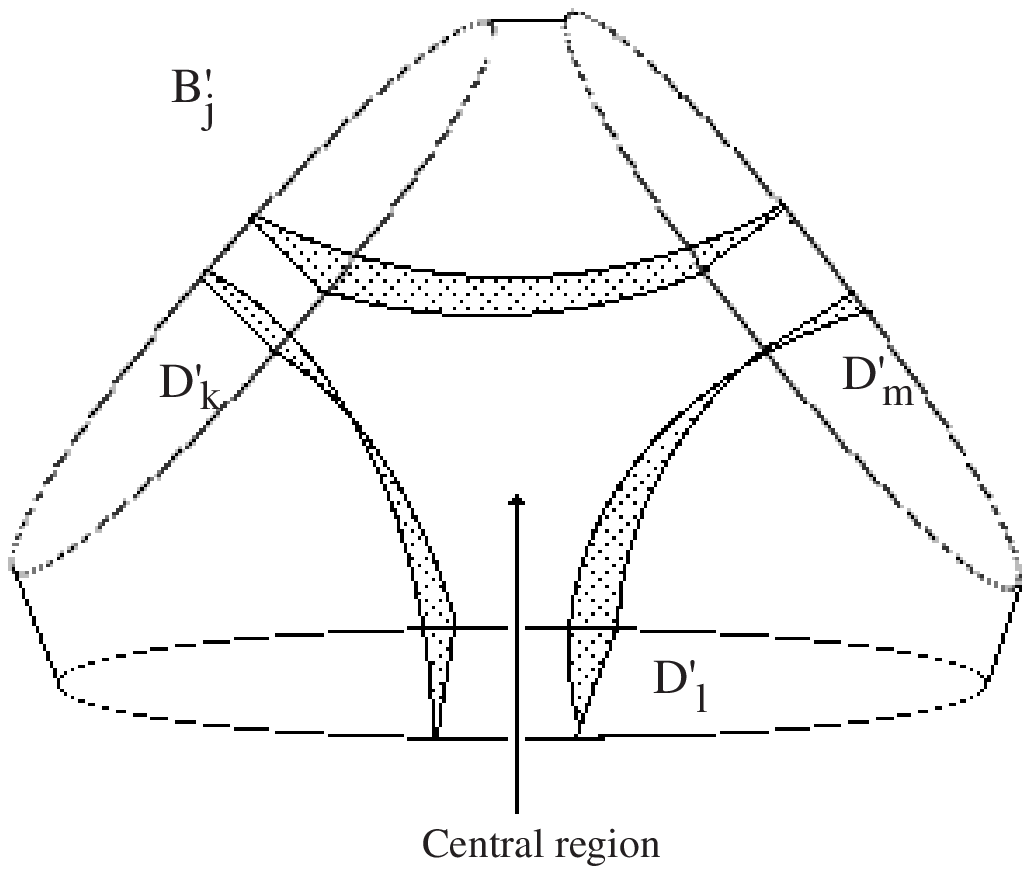}}}
\centerline{Figure 5.}

\vbox{{\epsfysize180pt\epsfbox{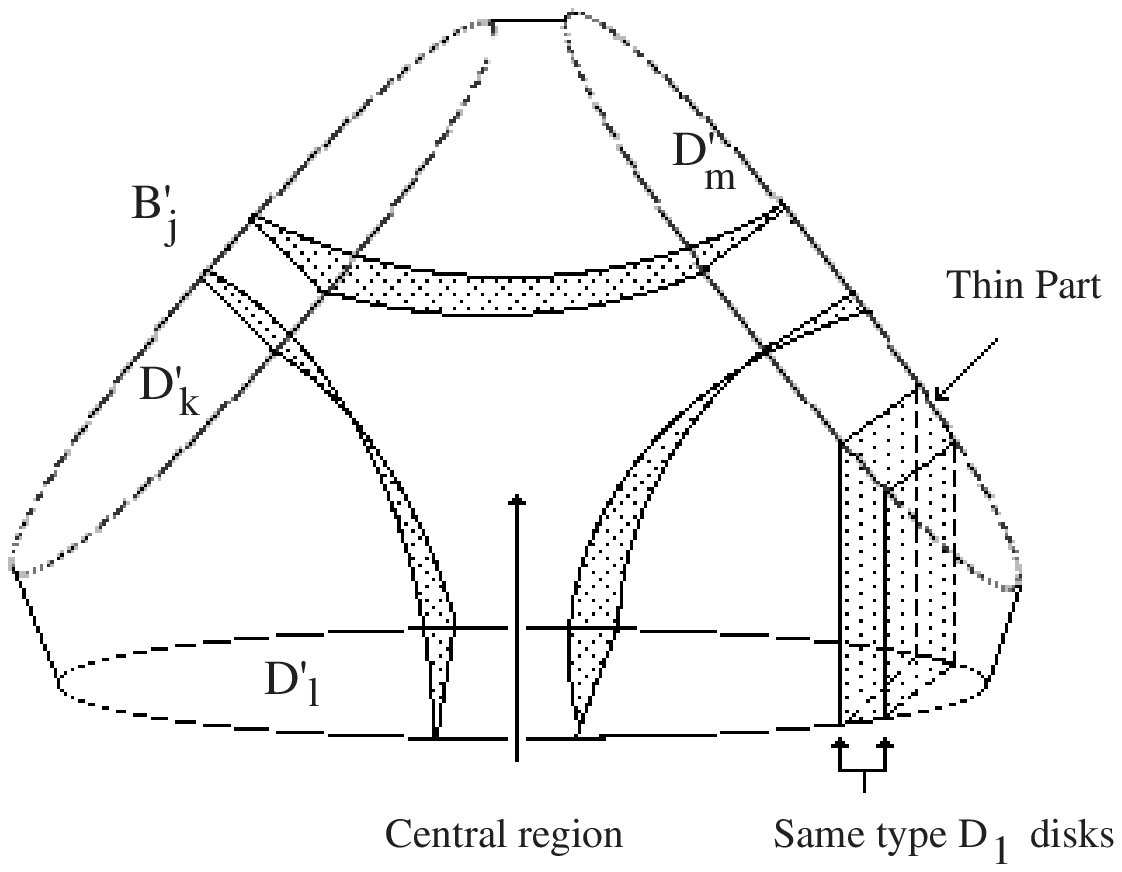}}}
\vskip10pt
\centerline{Figure 6. }

\vskip20pt

In Figure 7 below we display a schematic picture of a thick peripheral region. Note that in general they
can be more complicated.

\vbox{{\epsfysize260pt\epsfbox{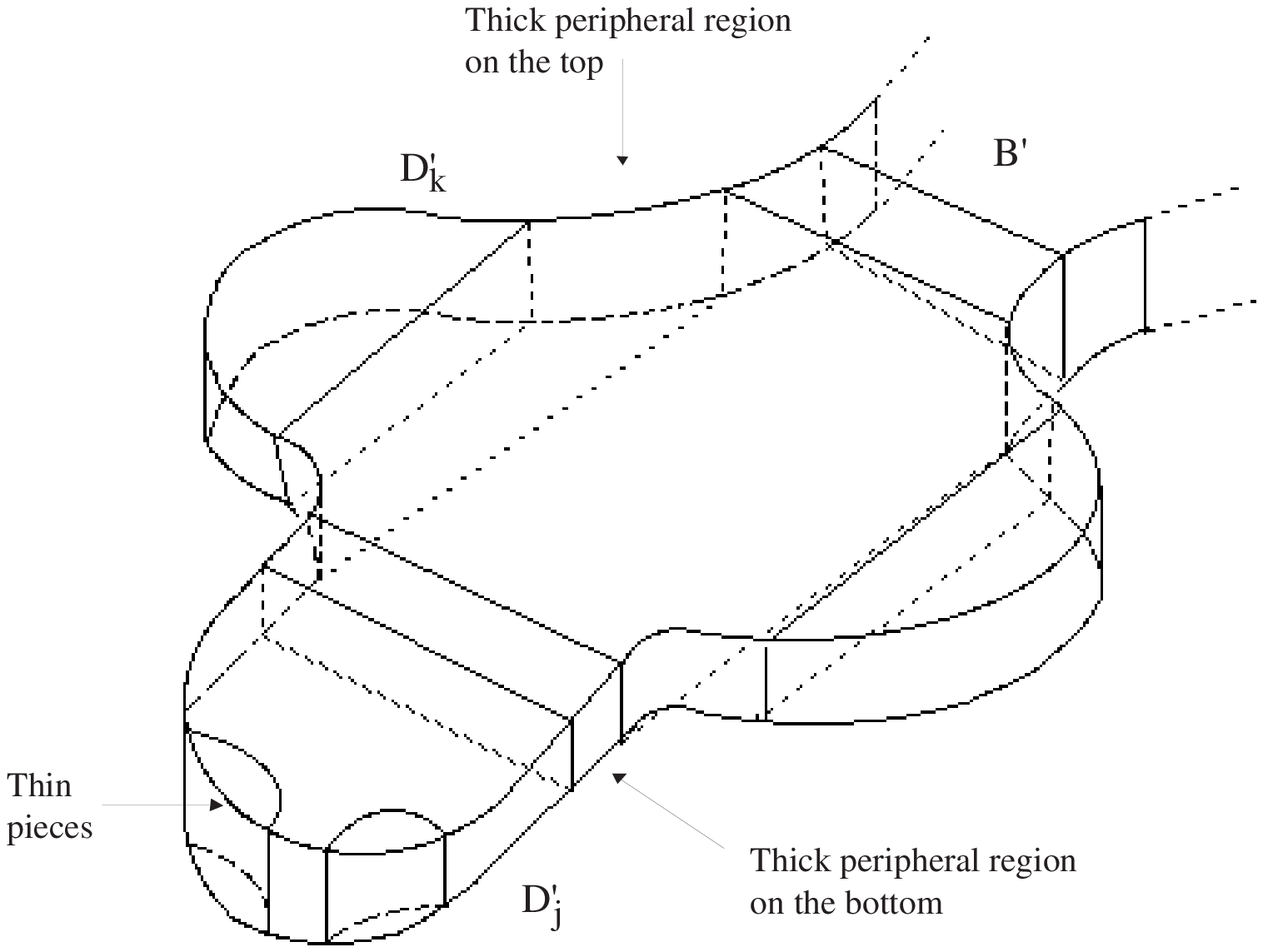}}}
\centerline{Figure 7. }

\vskip20pt

\begin{lemma}
\label{boundedarea}
There are finitely many possibilities for the central and the thick peripheral regions in $H_{1}$
which are completely determined by ${\cal D}_1$ and ${\cal D}_2$ only. In particular, their
number, as well as the volume and the diameter of any of them, are uniformly bounded above
by constants $N \geq 0$, $K \geq 0$ and $d \geq 0$ respectively.
\end{lemma}

\begin{proof}
We observed above that, in any solid pair of pants $B_{k}$, the complementary components
of the union of all stacks are precisely the  thick regions in $B_{k}$. Since the stacks are in
one to one correspondence with the disk types, the claim  follows directly from Proposition
\ref{finitetype}.

\end{proof}

\vskip8pt

\begin{lemma}
\label{uniquecenter}
(a)  Every solid pair of pants $B'_l$ has a unique central region.

\noindent(b)  For every disk $D'_j$ from ${\cal D}'_1$ which lies on the boundary of $B'_l$ the
intersection of $D'_j$  with the central region of $B'_l$ is connected.
\end{lemma}

\begin{proof}
(a) For every disk $D_i$ from ${\cal D}_1$ any connected component  $\Delta$ of $D_i \cap B'_l$ cuts
$B'_l$ into two distinct connected components. Hence, if $\Delta$ misses one of the three disks from
${\cal D}'_1$ which lie on the boundary of $B'_l$, say $D'_j$, then this disk $D'_j$ can intersect
only one of the two connected  components.

Now, note that by Definition \ref{regions} any two distinct thick regions in $B'_l$ are connected by a path
$\gamma$ which crosses at least one thin region, and hence, in the boundary of this thin region, $\gamma$
crosses a component $\Delta$ as above. This shows that at  most one of the two thick regions can be central.

To show the existence of a central region we first consider a connected component  $\Delta$ of
$D_i \cap B'_l$ which meets all three disks from ${\cal D}'_1$ that lie on the boundary of $B'_l$.
Such a $\Delta$ can not be contained in a thin or in a thick peripheral region, so that a central region
must exist. If there is no such $\Delta$, then, as shown above, each $\Delta$ cuts $B'_l$ into a connected
component that meets only two of the three boundary disks from ${\cal D}'_1$, and a second  connected
component that meets all three boundary disks. It follows directly that the intersection of these second connected
components, for all $\Delta$, is a single thick part which must meet all three boundary disks. Hence there
exists a central region to which this part belongs.

\noindent (b)  We observe that the subdisk $\Delta$ on the boundary of a thin region, as above, intersects a
disk $D'_j$ in at most one  arc. Hence we can apply the same arguments as in case (a) to any of the disks
$D'_j$ on the boundary of $B'_l$.

\end{proof}

\vskip8pt

A maximal connected union $P'$ of thin or thick peripheral parts of $B'_{l}$ is called {\em peripheral component}
of $B'_{l}$. Notice that any such peripheral component $P'$ meets precisely two disks $D'_{i}$ and $D'_{j}$ from
the collection ${\cal D}'_1$. It follows from the proof of Lemma \ref{uniquecenter} that the intersections
$P ' \cap D'_{i}$ and $P ' \cap D'_{j}$ are subdisks, and that $P'$ meets the closure of its complement $B'_{l} - P'$
in a subdisk $\Delta$ of some $D_{i}$ from ${\cal D}_1$, where $\Delta$ belongs to a thin part of $P'$. Hence the
boundary $\partial P'$ consists of $\Delta$, of $P' \cap D'_{i}$ and $P' \cap D'_{j}$, and of a band $A$ that has as
boundary two ``long" arcs $\alpha_{i} \subset \partial D'_{i}$, $\alpha_{j} \subset \partial D'_{j}$, and two ``short"
arcs $\beta, \beta' \subset \partial \Delta$.

\begin{lemma}
\label{peripheralvolume}

(a)  The arcs $\alpha_{i} = A \cap \partial D'_{i} $ and $\alpha_{j} = A \cap \partial D'_{j}$  meet exactly
the same sequence of disks from ${\cal D}_1$.

\noindent (b)  The number of disk pieces in the subdisks $P' \cap D'_{i}$ and $P' \cap D'_{j}$ is equal.

\end{lemma}

\begin{proof}
(a)  The band $A$ is topologically a disk (since $P'$ is a subball of the 3-ball $B'_{l}$), and we work
with the assumption that ${\cal D}_1$ and ${\cal D}'_1$ are tight, so that their boundary curves intersect
only essentially. Hence $\partial {\cal D}_1$ meets  $A$ in a collection of parallel arcs with one endpoint
on $\alpha_{i}$ and the other on $\alpha_{j}$.

\noindent
(b)  We observe that $P'$ may very well contain thick peripheral regions, so that the pattern of intersection
arcs on $P' \cap D'_{i}$ and on $P' \cap D'_{j}$ may be quite different. However, it follows directly
from (a) that the number of intersection arcs on $P' \cap D'_{i}$ and on $P' \cap D'_{j}$ must agree,
which implies the claim.

\end{proof}

\vskip10pt

Imagine the disk  $D'_j$ in a horizontal position, so that it  is part of the boundary of  an adjacent solid pair
of pants above it, and a second adjacent solid pair of pants below it. Both of these solid pairs of pants are
from the collection  ${\cal B}'$ defined  above. We call the intersection of $D'_j$ with the central region
from the top solid pair of pants the {\em top central subdisk}, and the one from the bottom the {\it bottom
central subdisk}. We measure the {\em distance} between them by counting the number of transverse
intersections with the disk system ${\cal D}_1$ of any path in $D'_j$ connecting the two central subdisks,
and taking the minimum over all such paths.

We define the {\em extended top central region} (and similarly the  {\em extended bottom central region})
to be the central region of the top solid pair of pants together with all parts from the bottom solid pair of
pants which are adjacent to the top central subdisk.

\vskip10pt

\begin{proposition}
\label{boundeddistance}
If ${\cal D}'_1$ and ${\cal D}'_2$ satisfy the double rectangle condition, then the distance between top
and bottom central subdisks on any of the disks $D'_j$ in ${\cal D}'_{1}$ is uniformly bounded from
above by a constant $c \geq 0$.

\end{proposition}

\begin{proof}
Since ${\cal D}'_1$ and ${\cal D}'_2$ satisfy the double rectangle condition, we can apply Lemma \ref{lemma1}
to show that every disk $E_i$ from the system ${\cal D}_2$ must intersect every adjacent triple from the system
${\cal D}'_1$ in some arc $\alpha_h \subset \partial E_i$. We consider in particular the four adjacent triples which are
contained in the union of the two solid pairs of pants $B'_l, B'_m$ adjacent to the disk $D'_j$ on the top and on
the bottom.

If the top central subdisk and the bottom central subdisk intersect in $D'_j$, then their distance is by
definition 0. In the case where the top and the bottom central subdisks of $D'_j$ are disjoint, we observe
that the extended top and bottom regions in $B'_l \cup B'_m$ are separated by pairs of parts, one on the
top, one on the bottom, which  belong to peripheral components of $B'_l$ and of $B'_m$. In particular,
the union of these pairs of parts meets only two of the four disks from ${\cal D}'_1$ which lie on the
boundary of $B'_{l} \cup B'_{m}$.

Hence for at least one of the above four adjacent triples, the corresponding  arc $\alpha_{h}$
intersects both, the top and bottom extended central regions. As a consequence, the distance
between the top and bottom central subdisks on $D'_j$ is bounded above by the minimal
number of intersections with ${\cal D}_1$ of any curve from ${\cal D}_2$. We will denote
this upper bound which depends only on  ${\cal D}_1$ and ${\cal D}_2$ by $c$.

\end{proof}

\vskip10pt

\section{Dual trees}
\label{dualtree}

\vskip15pt

For every disk $D'_j$ from ${\cal D}'_1$ we consider a graph whose vertices are in one to one
correspondence with the disk pieces of $D'_j$, and whose edges are in one to one correspondence
with the intersection arcs $\alpha_{i} \subset D'_j \cap {\cal D}_1$.  Each $\alpha_{i}$ cuts $D'_j$
into two distinct connected components. Hence the above graph is a tree, called {\it the dual tree}
$T'_{j}$.

We measure the {\it distance} in $T'_{j}$ by the usual simplicial metric, i.e. by associating to
every edge the length 1.  The {\em volume} of a subtree of $T'_{j}$ is given by the number
of vertices contained in the subtree. The {\em area} of a subdisk is the number of disk
pieces in the subdisk, which is equal to the volume of the corresponding dual subtree.

The top and bottom central subdisks of $D'_j$ define {\em top and bottom central subtrees} of $T'_{j}$.
The complementary components are called {\em top or bottom peripheral subtrees} of $T'_{j}$. Similarly,
any thick peripheral region in the adjacent top or bottom solid pair of pants defines, by way of intersection
with $D'_j$, {\em top or bottom thick peripheral subtrees} of $T'_{j}$.

\vskip10pt

\begin{remark}
\label{translate} \rm
Proposition \ref{boundeddistance} shows that the distance in $T'_{j}$ between the top and the
bottom central subtrees is uniformly bounded by the constant $c > 0$, if the double rectangle
condition is satisfied by ${\cal D}'_1$ and ${\cal D}'_2$.

\end{remark}

\vskip10pt

\begin{lemma}
\label{finitenessconditions}
For any real number $b > 0$ in any of the $T'_{j}$ the volume and the number of
complementary components of the $b$-neighborhood  of the top or of the bottom
central subtree, or of any thick peripheral subtree, are uniformly bounded above by some
constant $k = k(b) > 0$.

\end{lemma}

\begin{proof}
The valence of a given vertex in the tree is exactly the number of intersection arcs of the corresponding
disk piece with ${\cal D}_1$, which in turn is fixed for all the disk pieces from the same disk type. Since
there are only  finitely many disk types and finitely many thick parts, which are all determined
a priory  by ${\cal D}_1$ and ${\cal D}_2$, see  Proposition \ref{finitetype}, it follows directly
that there is an upper bound  $b_0$ on the valence of any vertex in $T'_{j}$. But then the
volume as well as the number of complementary components of the $b$-neighborhood of any finite
subtree is clearly bounded above by $b_0^b$ times the volume of that subtree, where the latter is bounded
uniformely in terms of the constant $K$ from Lemma \ref{boundedarea}.

\end{proof}

\vskip10pt

To continue the proof we need to define the following class of subtrees of any $T'_{j}$:

\vskip10pt

A subtree $R'_{j}$ of $T'_{j}$ will be called a {\em red subtree}, if it satisfies the following conditions:

\noindent (a) There is only one vertex, the {\em root} of $R'_{j}$, which is adjacent to some edge
contained in $T'_{j}$ but not in $R'_{j}$, and this edge is unique.  In other words, $R'_{j}$ is
obtained from $T'_{j}$ as connected component after removing a single edge.

\noindent (b) The subtree $R'_{j}$ is disjoint from the top or from the bottom central subtree of $T'_{j}$.
In the first case $R'_{j}$ is called a {\em top} red subtree, and in the second a {\em bottom} red subtree.

\vskip10pt

\noindent We now describe a two-tiered method, called the {\it disk pushing procedure}, of how to pass

\begin{itemize}

\item[ (I)] from a bottom red subtree in one of the $T'_{j}$ to a particular top red subtree in the same
$T'_{j}$, and

\item[{\rm (II)}] from a top red subtree in $T'_{j}$ to a particular bottom subtree in an adjacent $T'_{k}$.

\end{itemize}

\noindent It is this procedure that allows us to uniformly bound the size of the thin parts of the disks in
${\cal D}'_1$ and thus it is a crucial tool for the proof of our main result.

\vskip10pt

\noindent (I).  Let $R'_{j}$ be a bottom red subtree of $T'_{j}$.  We define an {\em adjacent top red subtree}
$R''_{j}$ as follows: If $R'_{j}$ is disjoint from the top central subtree of $T'_{j}$, then we set $R''_{j} = R'_{j}$.
If the top central subtree intersects $R'_{j}$, then we consider the $(c + d)$-neighborhood $C \subset T'_{j}$
of the top central subtree, for $c$ as in Proposition \ref{boundeddistance} and Remark \ref{translate}, and $d$
denoting the maximal diameter of the top or bottom central subtree in any of the dual trees $T'_{j}$
(see Lemma \ref{boundedarea}). Note that by Lemma \ref{boundedarea} the bound $d$ depends only on ${\cal D}_1$
and ${\cal D}_2$.  In this case we define $R''_{j}$ to be the complementary component of $C$ in
$R'_{j}$ which has largest volume. If there is more than one maximal volume component, pick any of them at
random.  It is immediate that $R''_{j}$ is a top red subtree as defined above.

Notice that,  in case $R''_{j} \neq R'_{j}$, since the bottom central subtree lies outside $R'_{j}$, the ``old"
root vertex, the one of $R'_{j}$, must either be contained in the top central subtree, or in the path in the
tree $T'_{j}$ which connects the bottom to the top central subtree. Thus Proposition \ref{boundeddistance}
implies that the root vertex of $R'_{j}$ is contained in $C$.  In particular we obtain the crucial fact that
$R''_{j}$ agrees with one of the complementary components of $C$ in $T'_{j}$, and not just in $R'_{j}$.

\vskip8pt
\noindent
(II).  If $R''_{j}$ is a top red subtree, then it is disjoint from the central region of the solid pair of pants
$B'_{l}$ which is adjacent to $D'_j$ from the top. Hence it is contained in a peripheral component of
$B'_{l}$.  Thus, among the three boundary disks of $B'_{l}$ from the system ${\cal D}'_1$, there is
precisely one, say $D'_k$, which differs from $D'_j$, but meets the same parts from $B'_{l}$ as $R''_{j}$.
From Lemma \ref{peripheralvolume} we know that the corresponding boundary arcs of $D'_j$ and $D'_k$
cross  exactly the same sequence of disks from ${\cal D}_1$.

We consider the dual tree $T'_{k}$ for $D'_k$, and the subtree $\hat T'_{k}$ of $T'_{k}$ which meets
the same parts of $B'_{l}$ as $R''_{j}$.  If the root vertex of $R''_{j}$ is contained in a thin part, then
we define the {\em the subsequent bottom red subtree} of $R''_j$ by $R'_{k} = \hat T'_{k}$.  Note that
in this case the trees $R''_{j}$ and $R'_{k}$ may be different (due to the presence of thick peripheral
regions), but, by Lemma \ref{peripheralvolume}, their volumes must agree.

If the root vertex of $R''_{j}$ is contained in a thick peripheral region, then we define the subsequent
bottom red subtree $R'_{k} $ to be a maximal complementary component in $\hat T'_{k}$ of this thick
peripheral region.  In either case the tree $T'_{k}$ is easily seen to satisfy again the properties of a bottom
red subtree (see Fig. 8).

\vbox{{\epsfysize260pt\epsfbox{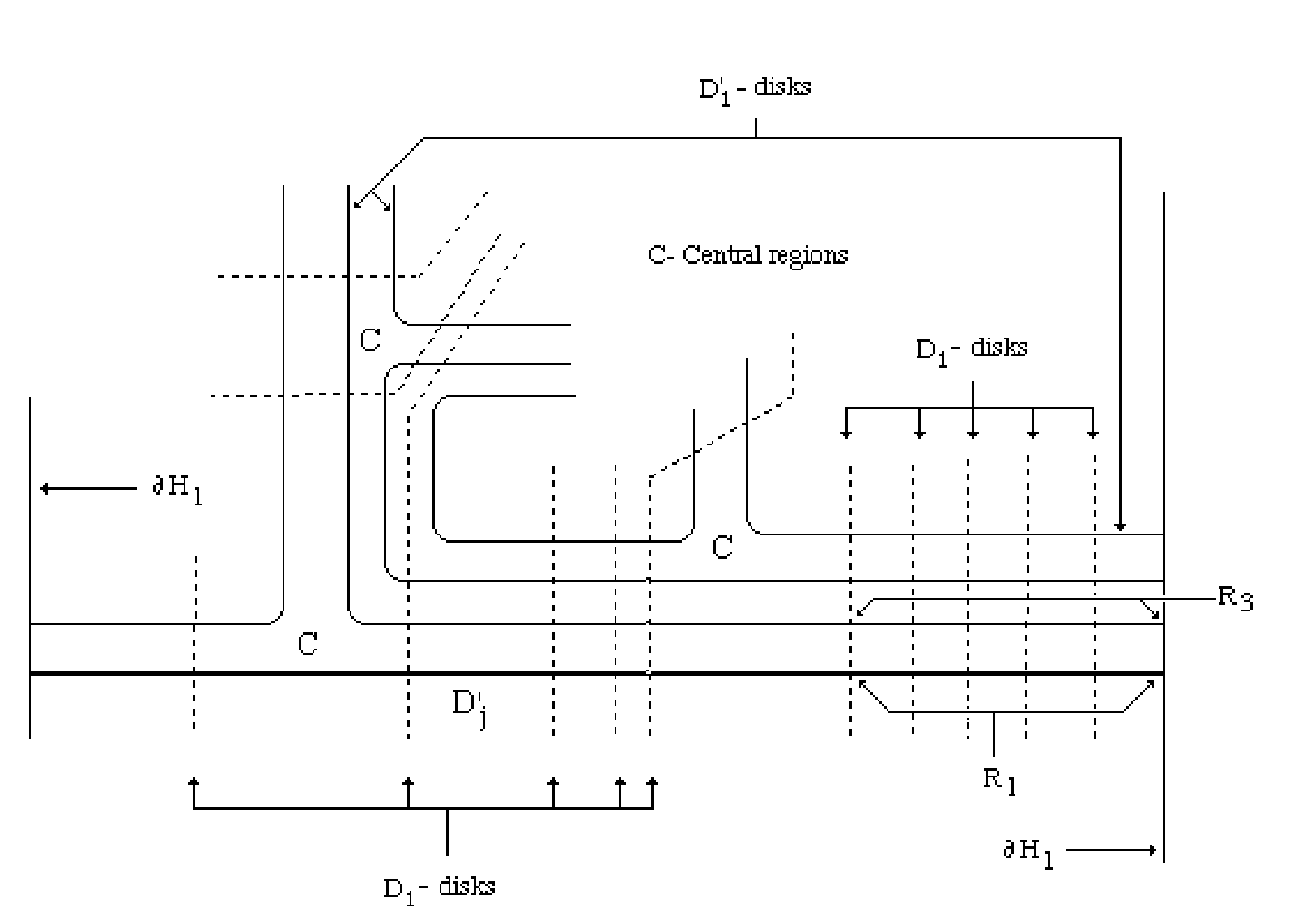}}}
\vskip30pt
\centerline {Figure 8.}

\vskip40pt

\begin{proposition}
\label{boundedseize}
For each pair of decomposing disk systems ${\cal D}_1 \subset H_{1}$ and ${\cal D}_2 \subset H_{2}$
there is an upper bound $a > 0$  such that, for any second pair ${\cal D}'_1 \subset H_{1}$ and
${\cal D}'_2 \subset H_{2}$ which satisfy the doubled rectangle condition, the area of any disk from
${\cal D}'_1$ with respect to ${\cal D}_1$ is bounded above by $a$.

\end{proposition}

\begin{proof}
By Lemma \ref{boundedarea} the systems ${\cal D}_1$ and ${\cal D}_2$ determine finitely many
possibilities for the thick regions in any solid pair of pants $B'_{l}$, and hence in particular for the
central subtrees for any of the adjacent disks $D_{j}$ from ${\cal D}'_1$. But, since any peripheral
subtree in the corresponding dual tree $T_{j}$ is a red subtree as defined above, our claim will be
proved if we show that the volume of any red subtree $R'_{j} \subset T'_{j}$ is bounded in terms
of ${\cal D}_1$ and ${\cal D}_2$.

Using the disk pushing procedure above we iteratively define a sequence of red subtrees $R_{n}$, starting
with $R_{1} = R'_{j}$, as follows:

\begin{itemize}

\item[(i)] If $R_{n}$ is a bottom red subtree, then $R_{n+1}$ is the adjacent top red subtree.

\item[(ii)]If $R_{n}$ is a top red subtree, then $R_{n+1}$ is the subsequent bottom red subtree.

\end{itemize}

Consider the sequence of volumes $r_{n}$ of the red subtrees $R_{n}$. This sequence is monotonically
decreasing (not necessarily strictly) for increasing $n$.  This follows directly from the definition of
$R_{n + 1}$ from $R_{n}$ by the disk pushing procedure.  In particular if $r_{n} = r_{n + 1} = r_{n + 2}$,
then the roots of the corresponding trees $R_{n}$, $R_{n + 2}$ are vertices in the corresponding dual trees
with the following property: The corresponding disk pieces (subdisks from the collection ${\cal D}'_1$)
do not belong to neither

\noindent (a) a thick peripheral region in both the top and bottom adjacent solid pairs of pants (from the
system  ${\cal B}'$), nor

\noindent (b) to the $(c + d)$-neighborhood of the central region of the top adjacent solid pair of pants.

As a consequence for any stationary subsequence $r_{n}, r_{n + 1}, r_{n + 2}, \dots, r_{n + 2k}$,
all root vertices of the corresponding trees $R_{n}, R_{n + 1}, R_{n + 2}, \dots, R_{n + 2k}$ belong to
distinct disk pieces $\Delta_{n}, \Delta_{n + 1}, \Delta_{n + 2}, \dots, \Delta_{n + 2k}$ which lie in the
stack of parallel disk pieces defined by a fixed disk type.  As all such stacks are finite (though not uniformly
bounded by ${\cal D}_1, {\cal D}_2$), it follows that any such stationary subsequence must be finite.

On the other hand, any time the value $r_{n+1}$ is strictly smaller than $r_{n}$, then the disk pushing
procedure for deriving $R_{n+1}$ from $R_{n}$ guarantees that $R_{n+1}$ coincides with a
complementary component in some of the $T'_{j}$ of either one of the bottom thick peripheral subtrees,
or of the $(c + d)$-neighborhood of one of the central subtrees.  The maximal number of such complementary
components is bounded above, by Lemma \ref{finitenessconditions}, by some $k = k(c + d)$, which only
depends on ${\cal D}_1$ and ${\cal D}_2$. Hence the number of values of the decreasing sequence of
areas $r_{n}$ is uniformly bounded.

It remains to observe that the quotient between two distinct values $r_{n}$ and $r_{n+1}$ is bounded
above in terms of ${\cal D}_1$ and ${\cal D}_2$ only: In fact, since in the definition of the adjacent
top red subtree, or of the subsequent bottom red subtree, we always chose a complementary subtree of
maximal volume, the inequality $\frac{r_{n} - k}{r_{n+1}} \leq k$ is valid for the value $k$ specified
above by Lemma \ref{finitenessconditions}. Recall here from Lemma \ref{finitenessconditions} that $k$
also bounds the volume of any thick peripheral or of the $(c + d)$-neighborhood of any central subtree in 
any of the $T'_{j}$.

This shows that the volume of any red subtree $R'_{j}$ is uniformly bounded above.

\end{proof}

From Proposition \ref{boundedseize} we immediately obtain  a proof of our main result Theorem \ref{mainthm} 
as stated in the Introduction. Notice that our proof is actually constructive, in that it describes a finite proceedure
which computes all complete decomposing systems which satisfy the double rectangle condition.

\begin{proof} We pick an arbitrary pair of complete	decomposing systems of disks ${\cal D}_1 \subset H_{1}$
and  ${\cal D}_2 \subset H_{2}$. By Proposition \ref{finitetype} there are finitely many disk types
with respect to these decomposing systems, which we can easily compute from the intersection pattern of
${\cal D}_1$ and ${\cal D}_2$ (see section \ref{disktypes}). By Lemma \ref{boundedarea} there are
finitely many possibilities for the thick regions, with an upper bound $N$ that only depends on
the already computed finite set of disk types.

We compute the upper bound $d$ for their diameter, the maximal length $c$ for any curve from ${\cal D}_2$, 
and the bound $k = k(c + d)$ as specified in the last proof. Then the formula $\frac{r_{n} - k}{r_{n+1}} \leq k$ 
from the last proof gives us the possibility to compute the largest possible area of any disk of the system ${\cal D}'_1$,
as in the decreasing sequence $r_{1}, r_{2}, \ldots$ the number of distinct values is bounded above by $k N$.

By symmetry we obtain a similar bound for the area of the disks from ${\cal D}'_2$, so that
there is only a finite number of candidates for  these systems, which can be directly computed from the
arbitrary chosen systems ${\cal D}_1$ and  ${\cal D}_2$.

\end{proof}

\begin{remark}
\label{yettobewritten}
\rm As mentioned in the Introduction it can be shown that the Casson-Gordon rectangle condition is 
generic in a precise meaning that uses Thurston's measure on the boundary of Teichm\"uller space. 
The analogous statement for the double rectangle condition, introduced in this paper, is not so clear.  
This is because one can define and impose an {\em anti double rectangle condition} as follows: 
The adjacent disk pairs from one side do not meet all four adjacent disk triples in a double pair of 
pants from the other side, but only three of them, and in place of the fourth one there is a repetition 
of one of the earlier triples, namely the one which is non-adjacent.  It is clear that the two conditions 
cannot be satisfied simultaneously. This anti double rectangle condition seems to be just as (non-)generic 
as the double rectangle condition.

A possible way to circumvent this difficulty is to consider the genericity of the set of
systems ${\cal D}_1, {\cal D}_2$ which (a) satisfy the Casson-Gordon rectangle condition,
and (b) have the property that ${\cal D}_1$ can be modified into a ``better" system 
${\cal D}'_1$ so that ${\cal D}'_1, {\cal D}_2$ satisfy the double rectangle condition.

An alternative approach,  which has implications
into other directions as well, is outlined as follows:

The role of the double rectangle condition is only to give an upper bound $c \geq 0$ on the maximal
distance  $c({\cal D}'_1, {\cal D}'_2)$ between the two central regions of any two adjacent pairs of
pants (compare Proposition \ref{boundeddistance}).  If we replace the double rectangle condition by
directly imposing such an upper bound on $c({\cal D}'_1, {\cal D}'_2)$ (defined in proper terms, so 
that the hypothesis becomes independent of the reference systems ${\cal D}_1, {\cal D}_2$ which 
are used to measure the quantity $c({\cal D}'_1, {\cal D}'_2)$), then the finiteness conclusion in our 
main Theorem \ref{mainthm} remains correct, and the proof stays virtually the same.  In this way 
we can define (despite Example \ref{example1}) for every Heegaard splitting which satisfies the 
Casson-Gordon rectangle condition for some disk systems ${\cal D}'_1, {\cal D}'_2$,  finitely many 
``preferred" such systems, namely those which have $c({\cal D}'_1, {\cal D}'_2)$ smaller than a given 
(sufficiently large) bound $c $.

\end{remark}

\eject

\section {References}

\vskip20pt

\noindent[CG] \hskip25pt A. Casson, C. Gordon;
{\it Reducing Heegaard
splittings}, Topology

\noindent  \hskip54pt  and its Applications {\bf 27}
(1987), 275 -- 283.

\vskip10pt

\noindent[He] \hskip28pt J. Hempel; {\it 3-manifolds
as viewed from the
curve complex},

\hskip33pt Topology {\bf 40} (2001), 631 -- 657.

\vskip10pt

\noindent[JR] \hskip28pt W. Jaco, H. Rubinstein; private
communication.

\vskip10pt

\noindent
[Sch] \hskip27pt S. Schleimer; {\it The disjoint
curve property},
preprint.

\vskip45pt

\obeyspaces   Martin Lustig \hskip149pt Yoav Moriah

\obeyspaces   Math\'ematiques (L.A.T.P.)
\hskip84pt Department of Mathematics

 \obeyspaces  Universit\'e d'Aix-Marseille III
 \hskip 66pt Technion, Haifa  32000,

\obeyspaces  Ave. E. Normandie-Niemen
\hskip 80pt Israel

\obeyspaces  13397 Marseille  20 \hskip 190pt

\obeyspaces  France \hskip 190pt

\medskip

\obeyspaces martin.lustig@univ.u-3mrs.fr \hskip75pt
ymoriah@tx.technion.ac.il

\end{document}